\documentclass[letterpaper, 10 pt, conference, twocolumn]{ieeeconf}  
\IEEEoverridecommandlockouts                            
\usepackage{cite}
\usepackage{amsmath,amssymb,amsfonts}
\usepackage[mathscr]{eucal}
\usepackage{algorithm}
\usepackage{algorithmic}
\usepackage{inputenc}
\usepackage{graphicx}
\usepackage{amsmath,tikz}
\usetikzlibrary{matrix}
\usepackage{verbatim}
\usepackage{hyperref}
\usepackage{textcomp}
\usepackage{subcaption}
\usepackage{capt-of}
\usepackage{subcaption}
\usepackage{booktabs}
\usepackage[font=small,labelfont=bf]{caption} 
\usepackage{lipsum}
\def\BibTeX{{\rm B\kern-.05em{\sc i\kern-.025em b}\kern-.08em
    T\kern-.1667em\lower.7ex\hbox{E}\kern-.125emX}}

\newtheorem{propo}{Proposition}[section]
\newtheorem{rema}{Remark}[section]
\floatname{algorithm}{Algorithm}
\DeclareMathOperator*{\argmin}{arg\,min}
\DeclareMathOperator{\tr}{tr}

\definecolor{Royalblue}{cmyk}{1,0.30,0.2,0.2}

\newcommand{\mz}{\color{black}}

\title{\huge{\textbf{Learning AR factor models}}}

\author{Francesca Crescente, Lucia Falconi, Federica Rozzi, Augusto Ferrante and Mattia Zorzi\thanks{F. Crescente, L. Falconi, F. Rozzi, A. Ferrante and M. Zorzi are with the Department of Information Engineering, University of Padova, Via Gradenigo 6/B, 35131 Padova, Italy. Emails: {\tt\small francesca.crescente@studenti.unipd.it}, {\tt\small lucia.falconi@studenti.unipd.it}, {\tt\small federica.rozzi@studenti.unipd.it}, {\tt\small augusto@dei.unipd.it}, {\tt\small zorzimat@dei.unipd.it}}}
 
\begin{document}

\maketitle
\thispagestyle{empty}
\pagestyle{empty}

\begin{abstract} 
We face the factor analysis problem using a particular class of auto-regressive processes. We propose an approximate moment matching approach to estimate the number of factors as well as the parameters of the model. This algorithm alternates a step of factor analysis and a step of AR dynamics estimation. Some simulation studies show the effectiveness of the proposed estimator.
 \end{abstract}

\section{Introduction} \label{sec:Introduction}
Factor models are among the first instances where statistical tools
have proved their power in providing sensible representations of a data
collection: the first contributions in this field date back to more than a century ago, \cite{spearman_1904,BURT_1909}.
In its simplest form the problem, aimed at extracting statistical commonalities
in multivariate data, may be reformulated as that of decomposing a positive definite covariance matrix $\Sigma$ as 
the sum 
\begin{equation}
\label{eq:basic-fa-problem}
   \Sigma=L+D
\end{equation}
where both $L$ and $D$ are positive semidefinite, $D$ is a diagonal matrix, and
$L$ has the lowest possibile rank. It turns out that, in general, this is a formidable problem for which
a rich stream of literature has been produced. We refer the reader to the
recent papers{\mz \cite{ning2015linear,bertsimas2017certifiably,8447502,8264253,bai2008large} where different principles for finding the decomposition in (\ref{eq:basic-fa-problem}) have been proposed and the papers \cite{Bottegal-Picci,DEISTLER_2015} where generalized formulations of the problem have been considered.} From the basic problem \eqref{eq:basic-fa-problem}, countless variations have been considered and studied. In particular, considerable effort has been devoted
to the dynamic case, see \cite{Anderson-Deistler-84,ANDERSON1985709,PICCI_PINZONI_1986,DEISTLER_1997,GEWEKE_DYNAMIC} and the references therein.

To practically compute a decomposition of type \eqref{eq:basic-fa-problem}
where the rank of $L$ is small, the typical strategy is to minimize a {\em proxy}
of the rank of $L$, i.e. the trace norm of $L$, \cite{MIN_RANK_SHAPIRO_1982,FAZEL_MINIMUM_RANK_2004,FAZEL_MIN_RANK_APPLICATIONS_2002}.
Trace norm regularization is also used in the so called latent-variable auto-regressive (AR) graphical models, \cite{TAC19,ZorzSep2016,ZORZI2019108516,ciccone2018robust} where we learn the spectral density of the model such that its inverse admits a ``sparse plus low-rank decomposition''. It is worth noting that data enters in these estimators through an approximate moments matching,  in a similar spirit of \cite{enqvist2007approximative}. This is a wise way to use moments.  
Indeed, in practice they are estimated from data and thus an estimation error is inevitable and  must be taken into account.

The natural dynamic extension to (\ref{eq:basic-fa-problem})
is to consider $\Sigma$, $L$ and $D$ as spectral densities of stationary stochastic processes, \cite{deistler2007}. In \cite{MFA}, factor analysis for moving average processes has been considered. The proposed estimator, however, matches exactly the ``noisy'' moments. As a consequence, the estimated decomposition is good provided that the number of data points is sufficiently large. 

In this paper, leveraging on the results in \cite{ciccone2018alternating,8447502}, we consider the problem of identification of the parameters of an AR model driven by a white noise whose covariance matrix
admits  a decomposition of the form \eqref{eq:basic-fa-problem} with the rank of $L$ being much smaller than the dimension of $\Sigma$.
This is an interesting situation because it corresponds to the case when
independent observation noises affect each channel of the observed AR process while a small number of common factors account for the information shared among the observations. \textcolor{black}{Our attention to AR processes is motivated by the fact that they can approximate arbitrary well any purely non deterministic process as long as the order is sufficiently high.} Our contribution is to propose an approximate moments matching method for the identification of the parameters
of these AR factor models. This method 
is based on alternating a step of factor analysis \textcolor{black}{(solved by resorting to the minimum trace proxy \cite{ciccone2018alternating})} and a step of AR dynamics estimation by means of moments matching. While this method mostly hinges on heuristic arguments,
it provides accurate estimations in high-dimensional instances: some simulations are 
described at the end of the paper that indeed show the performances
of the method for AR process with $40$ and $100$ channels.
\newline \indent The rest of the paper is organized as follows: Section \ref{sec:Problem} describes the problem formulation. Section \ref{sec:Algorithm} introduces the proposed algorithm. Numerical simulations and results are presented in Section \ref{sec:Simul}. Finally, Section \ref{sec:Conclusions} concludes the paper.

{\em Notation:} In this section we summarize and describe both the syntax and the semantics that will be used in the sequel.
\\Given a matrix $M$, $M^\top$ denotes the transpose; $|M|$ and $\tr(M)$ denote its determinant and trace (for a square $M$), respectively. The symbol $\mathbf{Q}_m$ denotes the space of real symmetric matrices of size $m$. If $M \in \mathbf{Q}_m$ is positive definite or positive semi-definite, then we write $M \succ 0$ or $M \succeq 0$, respectively. Moreover, we denote by $\mathbf{D}_m$ the space of diagonal matrices of size $m$. The symbols $\Vert\cdot\Vert_F$ and $\Vert\cdot\Vert$ stand for the Frobenius norm and Euclidean norm, respectively. 
With $f(\cdot)$ we denote the probability density function of a given random variable. The shorthand notation $x\,\bot\, y$ means that the random vectors $x$ and $y$ are independent.
We deal with Gaussian multivariate processes defined over the integers $\mathbb Z$.

\section{Problem Formulation} 
\label{sec:Problem}
\textcolor{black} {Consider, with an abuse of notation \footnote{ \textcolor{black}{For the sake of simplicity, we mix both the time and the $z$-domain representation.} } , the auto-regressive factor model: 
\begin{equation}
   y(t) = a(z)^{-1}[W_L v(t) + W_D w(t)] 
   \label{eq:ARMA_FA}
\end{equation}  }
where 
\begin{equation}
    a(z)=\sum_{k=0}^p a_{k}z^{- k}, \quad  a_{k}\in \mathbb{R},
    \label{eq:tau_def}
\end{equation}
$W_L \in \mathbb{R}^{m \times r}$, $W_D  \in \mathbf{D}_m$ is diagonal, and $p$ is the order of the model. The processes $v=\{v(t),\; t\in \mathbb{Z} \}$ and $w=\{w(t),\; t\in \mathbb{Z} \}$ are normalized white Gaussian noises of dimension $r$ and $m$ respectively; moreover, for all $t_1, t_2$, $v(t_1) \perp w(t_2)$.
The aforementioned model has the following interpretation: $v$ is the process which describes the $r$ factors, with  $r \ll m$, not accessible to observation;  $ a^{-1} W_L $ is the factor loading transfer matrix and $a^{-1} W_L v(t)$ represents the latent variable.  $a^{-1} W_D w(t)$ is the idiosyncratic noise describing the independent noises affecting each channel. 

Notice that 
\textcolor{black}{
\begin{equation}\label{eq:defu}
u(t):=  a(z) y(t) =   W_L v(t) + W_D w(t)  , 
\end{equation} }
is white Gaussian noise with covariance matrix given by \textcolor{black}{
$ \Sigma := W_L W_L^\top + W_D W_D^\top = L + D, $}
where $L:= W_L W_L^\top$ and $D:= W_D W_D^\top$.
We make the reasonable assumption that $\Sigma \succ 0$ (\textcolor{black}{since in most practical cases the noise affects the model in all the directions}) so that there exists $L_{\Sigma} \succ 0$ such that $ \Sigma = L_{\Sigma}L_{\Sigma}^\top.$
Then, (\ref{eq:defu}) may be written as
\textcolor{black}{
\begin{equation} 
    u(t) = a(z) y(t) =   L_{\Sigma} e(t)
    \label{eq:Taux} 
\end{equation}  }where $e=\{\,e(t),\; t\in\mathbb Z\,\}$ is an $m$-dimensional normalized white noise.
It follows that
\begin{equation} 
\textcolor{black}
{\bar{y}(t):= L_{\Sigma}^{-1}y(t) = a^{-1}(z) e(t)}
\label{eq:ARprocess}
\end{equation}
is still an AR process of order $p$. The process $\bar{y}(t)$ is obtained by stacking together 
the output of $m$ identical (scalar) filters driven by independent (scalar) white noises. Therefore $\bar{y}(t)$ may be viewed as a multivariate process with $m$ independent channels $\bar{y}_i(t)$ all of which feature the same probability description. This will be a key feature in what follows.

Assume now  to collect a finite length realization of $y$, say $\mathrm y^N=\{\, \mathrm{y}_1, \ldots ,\mathrm {y}_N\, \}$. Our aim is to estimate the corresponding factor model \eqref{eq:ARMA_FA} as well as the number of factors $r$.
The idea is to iteratively estimate $L$, $D$ and  $a(z)$ by pre-processing $\mathrm{y}^N$ through  $a$ and  $L_{\Sigma}^{-1}$,  respectively.
Finally, it is crucial to observe that there is an identifiability issue in this problem. Indeed, if we multiply $a(z)$ by an arbitrary non-zero real number $k$ and $L$ and $D$ by $k^2$, the model remains the same.
We can easily eliminate this uninteresting degree of freedom by normalizing the polynomial $a(z)$ so that from now on
we assume that $a_0=1$.



\section{Problem's solution}
\label{sec:Algorithm}

Our solution approach is based on an iterative algorithm that recursively estimates  $a(z)$, $L$, $D$ and $r$, until a certain tolerance is achieved. To easily explain our method, we firstly suppose $r$ is fixed. With such hypothesis, the proposed solution is presented in Algorithm \ref{algo:ARMAFA}. It receives as input the data $\mathrm{y}^N$, the order $p$ of $a(z)$, the number $r$ of factors, and the error tolerance $\varepsilon$. Given these inputs, it alternates the two following steps: 
\begin{enumerate} 
    \item the \textit{static factor analysis}, estimating the matrices $L$ and $D$;
    \item the \textit{AR dynamics estimation}, estimating the
  vector $\mathrm{a}:=[a_1 \ a_2 \dots a_p]^\top$ of the $p$ parameters of the polynomial $a(z)$ (recall that we have set $a_0=1$).
\end{enumerate}
These quantities are updated until the  difference between  two  consecutive  estimated  values of $L$, $D$ and $\mathrm{a}$ becomes smaller than a chosen threshold; more precisely, we impose that the mean square difference between the identification parameters in two successive steps is smaller than a given constant $\varepsilon$:
\begin{equation} \label{eq:stop_condition}
  \small{e:= \frac{\Vert L - L_{old} \Vert_F^2}{m^{2}}+
    \frac{\Vert D- D_{old} \Vert_F^2}{m}
    +\frac{\Vert \mathrm{a}- \mathrm{a}_{old}\Vert}{p} 
    \leq \varepsilon}
\end{equation}
where the subscript $old$ denotes the estimates of the previous iteration.
In addition, to ensure termination of the algorithm, we impose 
a maximum number $l_{max}$ of iterations. 

The two steps are explained hereafter, while the estimation of $r$ is addressed in Section \ref{sec:r}.

\begin{algorithm} 
\caption{Dynamic AR Factor Analysis}
\textbf{Input} $\mathrm{y}^N$, $p$, $r$, $\varepsilon$, $l_{max}$ \\
\textbf{Output}: $\mathrm{a}^{\circ}_r$,   $L^{\circ}_r$, $D^{\circ}_r$,
$\hat \Sigma_r$
\begin{algorithmic}[1] 
\STATE \textbf{Initialize} $\mathrm{a}$, ${L}$, ${D}$, $l=1$
\REPEAT
\STATE $\mathrm a_{old}=\mathrm a,\; L_{old}=L, \; D_{old}=D$
\STATE Data filtering $\mathrm u^N= \mathrm{a}\mathrm y^N$ \label{punto4}
\STATE Compute the sample covariance $\hat \Sigma $ from $\mathrm u^N$
\STATE \label{punto6} Compute $L,D$ s.t. $\mathrm{rank}(L)\leq r$, $D\in \mathbf D_m$ and\\ $\Vert \hat \Sigma-(L+D)\Vert_F$ is minimized 
\STATE Compute $L_{\Sigma}$ s.t.  $L + D = L_{\Sigma} L_{\Sigma}^{\top}$ 
\STATE \label{line8} Data ``disentangle'' $\bar{\mathrm y}^N= L_{\Sigma}^{-1}\mathrm y^N$
\STATE \label{punto9} Estimate  the AR coefficients  in $\mathrm{a}$  from $\bar{\mathrm y}^N$
\STATE Compute $e:= \frac{\Vert L - L_{old} \Vert_F^2}{m^{2}}+
    \frac{\Vert D- D_{old} \Vert_F^2}{m}
    +\frac{\Vert \mathrm{a}- \mathrm{a}_{old}\Vert}{p} $
\STATE $l=l+1$
\UNTIL  $e \leq \varepsilon $ \OR {$l \geq l_{max}$}
\STATE \quad $\mathrm{a}^{\circ}_r = \mathrm{a}$, $L^{\circ}_r = L$, $D^{\circ}_r=D $, 
$\hat \Sigma_r=\hat \Sigma$. 
\end{algorithmic}
\label{algo:ARMAFA}
\end{algorithm}

Some comments are in order.
\begin{enumerate}
\item
Formula $\mathrm u^N= \mathrm{a}\mathrm y^N$ in line 
\ref{punto4} of Algorithm \ref{algo:ARMAFA} has to be understood
as follows: consider the moving average filter $a(z):=1+ [z^{-1}\dots z^{-p}]\mathrm{a}$, $\mathrm u^N$ is the finite length trajectory obtained
by passing through the filter $a(z)$ the finite length trajectory $\mathrm y^N$ with zero initial conditions. Similarly for formula in 
line \ref{line8}: $\bar{\mathrm y}^N$ is the finite length trajectory obtained
by multiplying on the left side by $L_{\Sigma}^{-1}$ each vector of the finite length trajectory $\mathrm y^N$.
\item The aforementioned pre-processing steps are adaptive, indeed, at each iteration these operations changes according to the current $\mathrm a$ and $L_\Sigma$, respectively.
\item
Lines \ref{punto6} and  \ref{punto9} of Algorithm \ref{algo:ARMAFA}
correspond to other algorithms that are detailed in the following subsections.
\end{enumerate}

\subsection{Static factor analysis}
The static factor analysis problem is stated in \cite{ciccone2018alternating} as follows: for a given rank $r$ and a given matrix $\Sigma$ we want to find a positive semidefinite matrix $L$ with rank at most $r$ and a positive semidefinite diagonal matrix $D$ such that their sum is as close as possible to $\Sigma$. This can be formalized as:
\begin{equation}
    (L^*, D^*) := \argmin_{L \in {\mathcal{L}_{m,r}}, D \in \mathcal{D}_m} \Vert\Sigma - L - D \Vert ^2_F
    \label{eq:problem_static}
\end{equation}
where \textcolor{black}{$\mathcal{L}_{m,r} : = \{ X \in \mathbf{Q}_{m} : X \succeq 0, \; rank(X) \leq r\}$ and $\mathcal{D}_{m} : = \{ X \in \mathbf{D}_m : X \succeq 0 \}.$  }
To \textcolor{black}{efficiently} solve this problem we resort to Algorithm \ref{algo:alg1} that was first proposed and analyzed in \cite{ciccone2018alternating}. It receives as input the current matrix $\Sigma$ to be decomposed and the current value of the rank $r$, together with the error threshold $\varepsilon_{s}$. 
The estimation procedure is based on a coordinate descent type iterative algorithm. Such algorithm iterates between solving a minimization problem with respect to $L$ and a minimization problem with respect to $D$:
\textcolor{black}{$$
    L= \argmin_{L \in \mathcal{L}_{m,r}} \Vert \Sigma - L - D_{ old} \Vert^2_F, \;\;
    D = \argmin_{D \in \mathcal{D}_{m}} \Vert \Sigma - L - D \Vert^2_F
$$} 
where $D_{old}$ denotes the value of the diagonal matrix at the previous iteration. Notice that $P_{\mathcal{L}_{m,r}}$ and $P_{\mathcal{D}_{m}}$ in Algorithm \ref{algo:alg1} are the projector onto the sets $\mathcal{L}_{m,r}$ and $\mathcal{D}_{m}$ respectively. These projectors can be implemented very efficiently and robustly even for matrices $\Sigma$ with several hundreds of rows and columns.
\\The terminating condition is reached when $\Vert\Sigma - L-D\Vert ^2_F/\Vert\Sigma \Vert ^2_F \leq \varepsilon_{s}$ is satisfied.
\begin{algorithm} 	
		\caption{Static Factor Analysis}
	    \textbf{Input}: $\Sigma$, $r$, $\varepsilon_{s}$, $i_{\text{max}}$\\
	    \textbf{Output}: $L^*$, $D^*$
	    \begin{algorithmic}[1]
	    \STATE \textbf{Initialize}  $D$ randomly, $i=0$
		\WHILE{$\Vert\Sigma - L-D\Vert ^2_F/\Vert\Sigma \Vert ^2_F>\varepsilon_{s}$ \textbf{and} $i<i_{\text{max}}$}
		\STATE $L=P_{\mathcal{L}_{m,r}}(\Sigma-D)$
		\STATE $D=P_{\mathcal{D}_{m}}(\Sigma-L)$
		\STATE $i=i+1$
    	\ENDWHILE
    	\STATE $L^*=L, \; D^*=D$
    	\end{algorithmic}
    	\label{algo:alg1}
\end{algorithm}

\subsection{AR dynamics estimation}
\label{sec:AR}
The second step is the AR dynamics estimation. Given a finite-length realization $\bar{\mathrm{y}}^N=\{\bar{\mathrm{y}}(1)\ \dots\ \bar{\mathrm{y}}(N)\}$ of the AR process \eqref{eq:ARprocess}, the aim is to estimate the coefficients of the filter $a(z)$, namely $\mathrm{a} = [a_1 \ \cdots \ a_p]^\top \in \mathbb{R}^{p}$ (as we have fixed $a_0=1$).
To this aim we resort to the maximum-likelihood (ML)   principle and
compute the estimate as
\begin{equation}
    \mathrm{a}_{ML} := \argmin_{\mathrm{a}} \ell (\bar{\mathrm{y}}^N; \mathrm{a})
    \label{eq:ARproblem} 
\end{equation}
where the negative log-likelihood $\ell ( \bar{\mathrm{y}}, \mathrm{a})$ is defined as
\begin{equation*}\begin{split}
    \ell (\bar{\mathrm{y}}^N; \mathrm{a}) :=  - \log f \Big( \bar{\mathrm{y}}(N),.., \bar{\mathrm{y}}(p+1) | \bar{\mathrm{y}}(p),.., \bar{\mathrm{y}}(1) \Big).
\end{split}\end{equation*}
In other words, we estimate the parameters vector $\mathrm{a}$ in such a way that the  model 
\begin{equation}
\textcolor{black}{
\bar{y}(t)= a(z)^{-1}e(t)},
\label{eq:y_t}
\end{equation}
maximizes the likelihood of producing the finite trajectory $\bar{\mathrm{y}}^N$.

Firstly, we consider the scalar case $m = 1$. 
Since we are dealing with an AR model the solution can be obtained by standard arguments in closed form. In fact,
by taking \eqref{eq:tau_def} into account, we can rewrite \eqref{eq:y_t} as \textcolor{black}{$
\sum_{k=0}^{p}a_k \bar{y} (t-k)= e (t), $ }
so that \textcolor{black}{$ 
\bar{y} (t)=-\sum_{k=1}^{p} a_k \bar{y} (t-k)+ e(t).$}  Therefore,
\begin{equation}
f\bigl( \bar{y}(t) | \bar{y}(t-1)\dots \bar{y}(t-p) \bigr)  \sim  \mathcal{N} \bigg( -[\, 0\; \mathrm a^\top \,]Y(t),1 \bigg)\label{eq:prob}
\end{equation}
and
\begin{equation}
    \begin{split}
    f \bigl( \bar{y}(N)\dots \bar{y}(p+1)| \bar{y}(p)\dots \bar{y}(1) \bigr)=\hspace{28mm} \\
    \prod_{t=p+1}^N f \bigl( \bar{y}(t) | \bar{y}(t-1)\dots \bar{y} (t-p) \bigr)
    \end{split}
\label{eq:prob1}
\end{equation}
 where $Y(t)=[\, \bar y(t)\; \bar y(t-1)\dots y(t-p) \,]^\top$.
Then, the negative log-likelihood (up to constant terms) results
\begin{equation*}
\begin{split} 
\ell(\bar{\mathrm{y}}^N;\mathrm{a})  & = \frac{1}{2} \sum_{t=p+1}^{N} \left( [\, 1 \; \mathrm a^\top\, ] Y(t) \right)^2 \\
& =  \frac{1}{2} \sum_{t=p+1}^{N} \tr \bigg( [1\  \mathrm{a}^\top] Y(t) Y(t)^\top \begin{bmatrix} 1 \\ \mathrm{a} \end{bmatrix} \bigg) \\ 
& = \frac{1}{2}\ 
[1\  \mathrm{a}^\top] \left[\sum_{t=p+1}^{N} Y(t) Y(t)^\top\right] \begin{bmatrix} 1 \\ \mathrm{a} \end{bmatrix}
\end{split} 
\end{equation*}
We now define the matrix 
$$\hat{T} :=  \frac{1}{N-p}\sum_{t=p+1}^{N} Y(t) Y(t)^\top=
\begin{bmatrix} \tilde \tau_0 & z^\top \\z & \hat T_{22}\end{bmatrix}$$
where the partition is such that $\tilde \tau_0 $ is a scalar and
$z$ is a column vector.
In this way, we have 
 \begin{equation*} 
\ell ( \bar{\mathrm{y}}^N;\mathrm{a})   = 
\frac{N-p}{2}\ 
[ 1\  \mathrm{a}^\top] \hat{T}
\begin{bmatrix} 1 \\ \mathrm{a} \end{bmatrix}=\frac{1}{2}\left[ \tilde \tau_0 +2 z^\top\mathrm{a} + \mathrm{a}^\top \hat T_{22} \mathrm{a}\right]
\end{equation*}

Since $\ell ( \bar{\mathrm{y}}^N;\mathrm{a})$ is clearly convex in $\mathrm{a}$,
Problem \eqref{eq:ARproblem} is solved by annihilating the gradient of 
$\ell ( \bar{\mathrm{y}}^N;\mathrm{a})$ with respect to $\mathrm{a}$, i.e. by imposing
that
\begin{equation} \label{eq:der}
   \frac{\partial \ell ( \bar{\mathrm{y}}^N;\mathrm{a})}{\partial \mathrm{a}}=2 z^T +2 \hat T_{22}\mathrm a =0.
\end{equation}
Since $\ell ( \bar{\mathrm{y}}^N;\mathrm{a})$ is a quadratic form in  $\mathrm{a}$,
\eqref{eq:der} provides a closed form formula
which finally yields 
\begin{equation} 
\mathrm{a}_{ML}= - \hat T^{-1}_{22} z.
\label{eq:derzero3}
\end{equation}
Of course, the interesting case is the multivariate one i.e. $m>1$.
To address this case, we recall that the $m$ components $\bar{y}_i(t)$
of the vector process $\bar{y}(t)$ are independent scalar processes
i.e.  \textcolor{black}{$\bar{y}_k(t)\perp \bar{y}_l(s),\quad \forall k\ne l, \  \forall s,t,$}
and they all have the same probabilistic description i.e. all the $\bar{y}_k(t)$'s  have the same spectral density 
\begin{equation}\label{eq:specdencomp}
\Phi_{\bar{y}_k}(z) = \frac{1}{a(z)a(z^{-1})}, \quad \forall i=1,\dots, m.
\end{equation}

The multivariate case can therefore be addressed as that of a scalar process with 
$m$-times as many data. In fact, in view of \textcolor{black}{the independence of the components of $\bar{y}(t)$,}  the likelihood is
\begin{equation*}
\begin{split}
    f(&\bar{y}(N)\dots \bar{y}(p+1)|\bar{y}(p)\dots \bar{y}(1))=\hspace{32mm} \\
    & \hspace{18mm} \prod_{k=1}^m f\left(\bar{y}_k(N)\dots \bar{y}_k(p+1)|\bar{y}_k(p)\dots \bar{y}_k(1)\right).
    \end{split}
\end{equation*}
Moreover, in view of \eqref{eq:specdencomp},
we can repeat for each $k$ the argument that led to \eqref{eq:prob1}
to obtain an expression for
$f\left(\bar{y}_k(N)\dots \bar{y}_k(p+1)|\bar{y}_k(p)\dots \bar{y}_k(1)\right)$.
This yields
\begin{equation*}
\begin{split}
    f(&\bar{y}(N)\dots \bar{y}(p+1)|\bar{y}(p)\dots \bar{y}(1))=\hspace{32mm} \\
    & \hspace{18mm} \prod_{k=1}^m \prod_{t=p+1}^N f \bigl( \bar{y}_k(t) | \bar{y}_k(t-1)\dots \bar{y}_k (t-p) \bigr).
    \end{split}
\end{equation*}
We can  now repeat the previous computation and obtain
\begin{equation*} 
\ell ( \bar{\mathrm{y}}^N;\mathrm{a})  =\frac{1}{2} \sum_{k=1}^{m}  \ 
[1\  \mathrm{a}^\top] \left[\sum_{t=p+1}^{N} Y_k(t) Y_k(t)^\top\right] \begin{bmatrix} 1 \\ \mathrm{a} \end{bmatrix}
\end{equation*}
where $Y_k(t):=[\,  \bar{\mathrm{y}}_k(t)\;  \bar{\mathrm{y}}_k(t-1)\dots \bar{\mathrm{y}}_k(t-p)\,]^T$.
We now define the matrix 
{\small\begin{equation}
\label{eq:TeZfinali} 
\hat{T} := \frac{1}{m(N-p)}\sum_{k=1}^{m}\sum_{t=p+1}^{N} Y_k(t) Y_k(t)^\top=
\begin{bmatrix} \tilde \tau_0 & z^\top \\z & \hat T_{22}\end{bmatrix}
\end{equation}}where, as for the scalar case, the partition is such that $\tau $ is a scalar and
$z$ is a column vector.
In this way, we are exactly in the situation discussed for the scalar case and the solution is thus given again by
\eqref{eq:derzero3} with $T$ and $z$ now provided by \eqref{eq:TeZfinali}. Such a solution, however, is not guaranteed to correspond to a stable model (i.e. a model such that all the zeros of $a(z)$ are inside the unit circle). Notice that $\hat T$ is an estimate of the Toeplitz matrix $T=\mathbb E[Y_k(t)Y_k(t)^\top]\succ 0$. Although $\hat T\rightarrow T$ almost surely as $N\rightarrow \infty$, $\hat T\succ 0$ it is {\em not} Toeplitz for finite values of $N$. To address such an issue, we consider the biased estimate 
\begin{align}
\hat T_b=\left[\begin{array}{cccc}\tau_0 & \tau_1 & \ldots  & \tau_{p} \\ \tau_1 & \tau_0 & \ddots  & \\     &\ddots  &  \ddots & \tau_1\\ \tau_{p} & & \tau_1 &  \tau_0\end{array}\right]=\left[\begin{array}{cc}\tau_0 &  z_b^\top \\ z_b & \hat T_{b,22}\end{array}\right]
\end{align}
where \textcolor{black}{$\tau_l= \frac{1}{m N}\sum_{k=1}^m \sum_{t=k+1}^N \mathrm y_k(t) \mathrm y_k(t-l),\; l=0\dots p.$}
It is not difficult to see that $\hat T_b\succ 0$ generically. Accordingly, we can choose as estimate of $\mathrm a$:
\begin{align}\label{eq_a_ME}
\mathrm a_{ME}=-\hat T_{b,22}^{-1}z_b.
\end{align} It is worth noting that (\ref{eq_a_ME}) is the solution to a Yule-Walker equation \cite{stoica2005spectral}.
Accordingly, $a_{ME}(z)=1+[z^{-1}\dots z^{-p}]\mathrm a_{ME}$ is a stable polynomial. Hence, the estimated spectral density of each $\bar y_k(t)$ is $\Phi_{ME}(z)=(a_{ME}(z)a_{ME}(z^{-1}))^{-1}$.
\begin{propo} Let $\check \tau_l=\tau_l v^\top \hat T_b^{-1}v$ with $v=[\,1 \;0\dots 0\, ]^\top$. Then, $\Phi_{ME}$ is the unique solution to the following maximum entropy problem:
\begin{align}\label{pb_ME}\Phi_{ME}=&\underset{\Phi}{\mathrm{argmax}}\int_{-\pi}^\pi \log\det \Phi(e^{j\vartheta})\mathrm d \vartheta\\ & \hbox{ s.t. } \int_{-\pi}^\pi  e^{j\vartheta l}\Phi(e^{j\vartheta})\mathrm d \vartheta=\check \tau_l, \; l=0\dots p.\end{align}
\end{propo}
{\em Proof:} Let $\check T_b:=v^\top \hat T_b^{-1}v\hat T_b$. It is well known \textcolor{black}{(see for example \cite{stoica2005spectral})} that the solution to (\ref{pb_ME}) is $\Phi(z)=\sigma^2( a(z)a(z^{-1}))^{-1}$ with $a(z)=1+[z^{-1}\dots z^{-p}]\mathrm a$ such that
\begin{align}\label{YW_rescaled} \check T_b \left[\begin{array}{c}1 \\ \mathrm a\end{array}\right]=\left[\begin{array}{c}\sigma^2 \\ \mathrm 0\end{array}\right].\end{align} Notice that in (\ref{YW_rescaled}) is a system of $p+1$ equations. Consider the subsystem composed by the second equation up to the last equation: since $\check T_b$ 
is invertible, its solution is (\ref{eq_a_ME}). It remains to show that $\sigma^2=1$. Substituting (\ref{eq_a_ME}) in the first equation, we have 
\begin{align*} \sigma^2&=v^\top \check T_b  [\,1 \;\mathrm a_{ME}^\top\,]^\top=(\tau_0+z_b^\top \mathrm a_{ME})(v^\top \hat T_b^{-1}v)\\
& =(\tau_0-z_b^\top \hat T^{-1}_{b,22}z_b)(v^\top \hat T_b^{-1}v)=1
\end{align*}
where the last equality is due by the fact that $\tau_0-z_b^\top \hat T^{-1}_{b,22}z_b$ is the Schur complement of the block $\hat T_{b,22}$ of $\hat T_b$. $\; \; \blacksquare$

\begin{figure*}
   \centering
   \begin{subfigure}[b]{0.4\textwidth}
       \includegraphics[width=1.1\textwidth, height = 0.17\textheight]{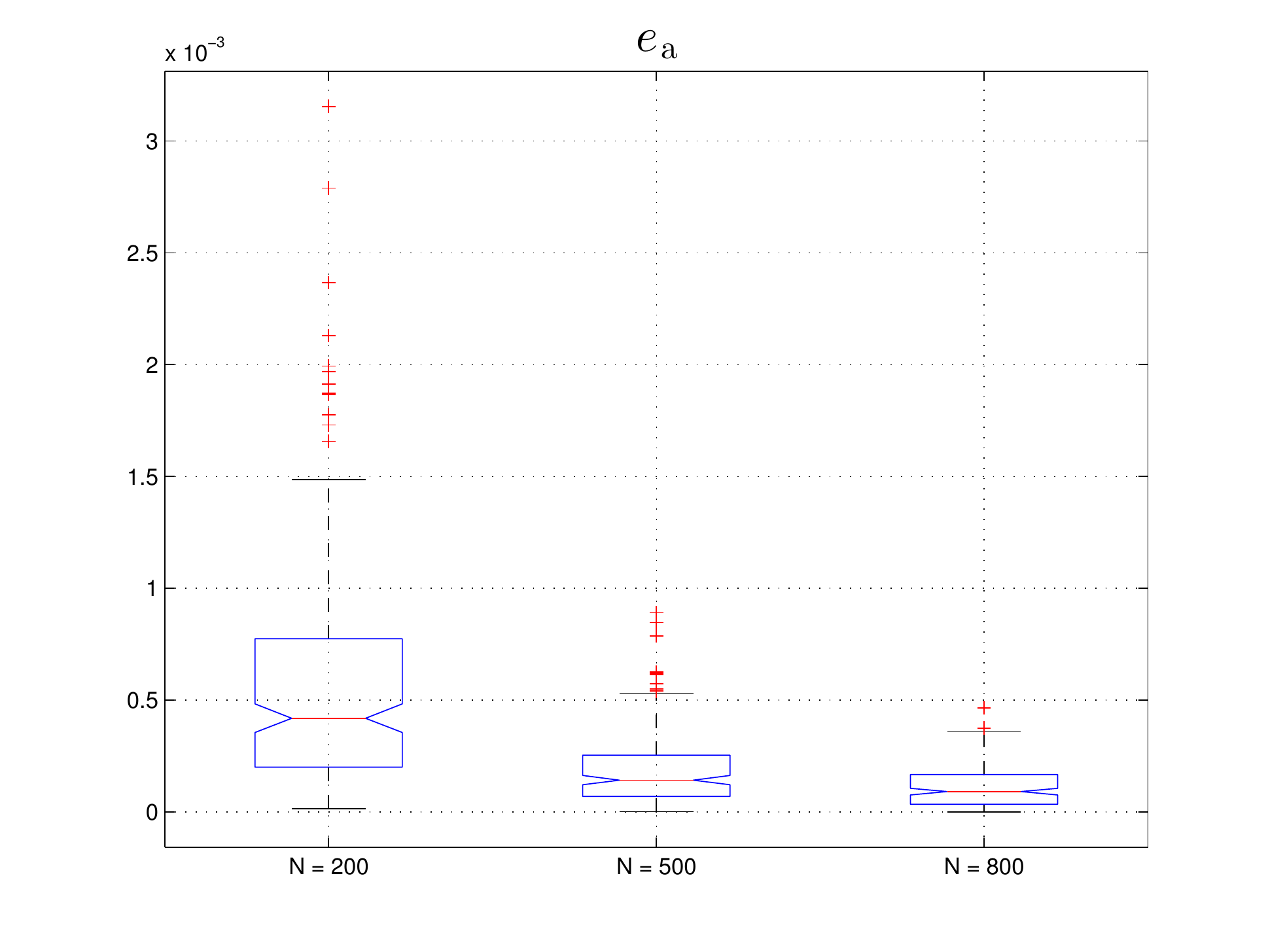}
       \caption{}
       \label{fig:m40_errTau}
   \end{subfigure}
   \begin{subfigure}[b]{0.4\textwidth}
       \includegraphics[width=1.1\textwidth, height = 0.17\textheight]{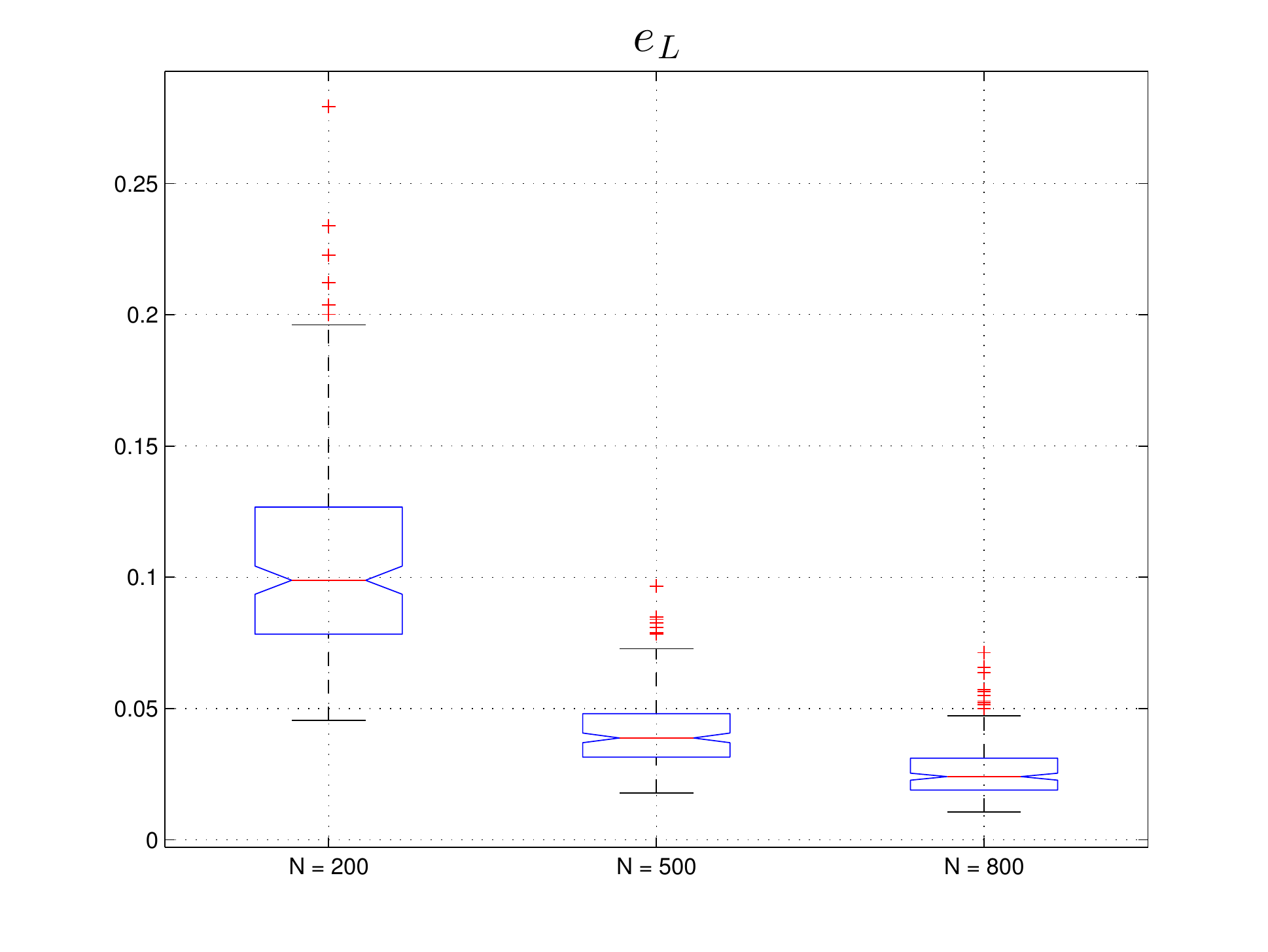}
       \caption{}
       \label{fig:m40_errL}
    \end{subfigure}
    \begin{subfigure}[b]{0.4\textwidth}
       \includegraphics[width=1.1\textwidth, height = 0.17\textheight]{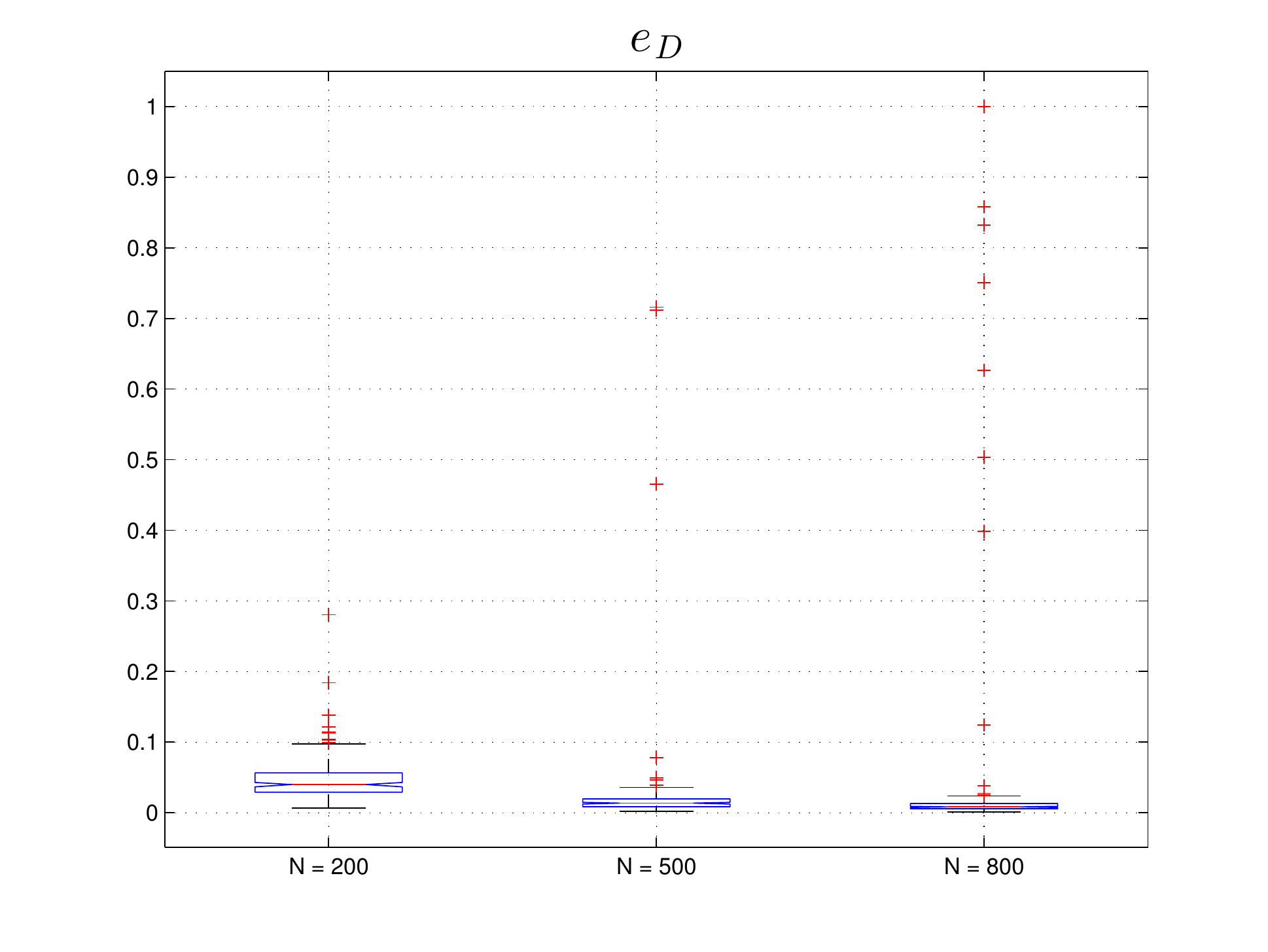}
       \caption{}
       \label{fig:m40_errD}
   \end{subfigure}
   \begin{subfigure}[b]{0.4\textwidth}
       \includegraphics[width=1.1\textwidth, height = 0.17\textheight]{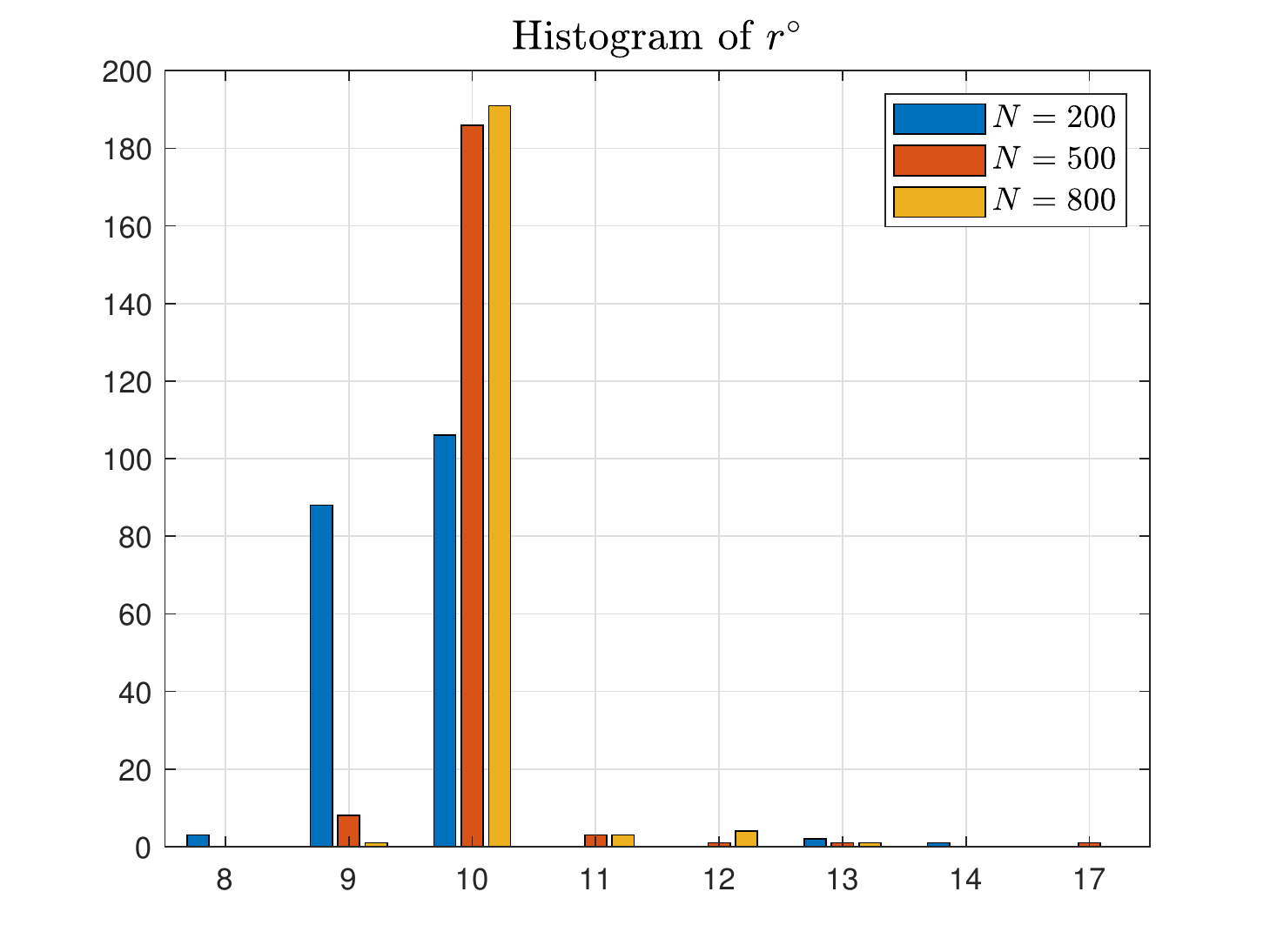}
       \caption{}
       \label{fig:m40_bar}
    \end{subfigure}
    \caption{Case $m=40$,  $r=10$ and $p=5$ with $N=200$, $N=500$ and $N=800$. Figure (a), (b) and (c) show the box-plot of the error $e_{\mathrm{a}}$, $e_{\mathrm{L}}$ and $e_{\mathrm{D}}$ respectively; (d) is the bar-plot of estimated $r^\circ$.}
    \label{fig:m40}
\end{figure*}

Algorithm \ref{algo:ARproc} summarizes the AR estimation procedure.
\begin{algorithm} 	\caption{AR dynamics Estimation}
        \textbf{Input}: $p,\; \bar{\mathrm{y}}^N$\\
        \textbf{Output}: $\mathrm{a}_{ME}$ 
        \begin{algorithmic}[1] 
        \STATE Let $Y_k(t):=[\bar{\mathrm{y}}_k(t)  \dots \bar{\mathrm{y}}_k(t-p)]^\top$, $k=1\dots m$.
        \STATE Compute $\hat T_b$ and thus $\hat T_{b,22}$, $z_b$ from $Y_k(t)$
\STATE $\mathrm{a}_{ME} =- \hat T_{b,22}^{-1} z_b$
 \end{algorithmic}
 \label{algo:ARproc}
\end{algorithm}

\begin{rema} It is clear that $\mathrm a^\circ_r$ in Algorithm \ref{algo:ARMAFA} matches the rescaled moments of $\bar{\mathrm y}^N=L_{\hat \Sigma_r}^{-1} \mathrm y^N$. Moreover, $L_r^\circ+D_r^\circ$ approximately matches the zeroth moment of $\mathrm u^N=\mathrm a^\circ_r \mathrm y^N$. Accordingly, Algorithm \ref{algo:ARMAFA} is an approximate moment matching estimator for an AR factor model of type (\ref{eq:ARMA_FA}).
\end{rema} 

\subsection{Estimation of $r$}
\label{sec:r}
As regards the estimation of the number $r$ of factors in \eqref{eq:ARMA_FA}, 
we propose the following procedure. We start from the most parsimonious model
with only a single factor and increase the number of factors until the difference 
between $\hat{\Sigma}_r$ and the sum $L^{\circ}_r+D^{\circ}_r$ 
(where $\hat{\Sigma}_r$, $L^{\circ}_r$ and $D^{\circ}_r$ are
outputs of Algorithm \eqref{algo:ARMAFA}) is sufficiently small to be 
explained by the estimation error (due to finiteness number of data) of the sample covariance $\hat{\Sigma}_r$. More precisely, starting from ${r}=1$, Algorithm \ref{algo:ARMAFA} is iteratively applied, increasing ${r}$ at each step. The stopping criterion is defined in the sequel. Let $\hat{\Sigma}_{{r}}$, $L_{{r}}^{\circ}$ and $D_{{r}}^{\circ}$ be the achieved values using ${r}$.
Consider the Kullback-Leibler divergence between $\hat{\Sigma}_{{r}}$ and  $\Sigma_{{r}}^{\circ} := L_{{r}}^{\circ} + D_{{r}}^{\circ}$, defined as:
\begin{equation*} \label{KL_divergence}
    \mathcal{D}_{KL} (\Sigma_{{r}}^{\circ} \Vert \hat{\Sigma}_{{r}}) := \frac{1}{2} ( -\log |\Sigma_{{r}}^{\circ}| + \log | \hat{\Sigma}_{{r}} | + \tr (\Sigma_{{r}}^{\circ} \hat{\Sigma}_{{r}}^{-1} ) - m).
\end{equation*}
Then, the value of ${r}$ is increased until $\mathcal{D}_{KL} (\Sigma_{{r}}^{\circ} \Vert \hat{\Sigma}_{{r}})$ becomes smaller than a given tolerance $\delta$, i.e.
\begin{equation}\label{eq:stop_criterion}
    \mathcal{D}_{KL} (\Sigma_{{r}}^{\circ} \Vert \hat{\Sigma}_{{r}}) \leq \delta
\end{equation}
 and such value of $r$, denoted by $r^\circ$, is the estimate of the rank. Clearly, the optimal model is given by $\mathrm a^\circ_{r^\circ}$, $L^\circ_{r^\circ}$ and $D^\circ_{r^\circ}$.

It is worth noting that $\mathcal{D}_{KL} (\Sigma_{{r}}^{\circ} \Vert \hat{\Sigma}_{{r}})$ measures how well the model with $r$ factors explains the data. Accordingly, the stop criterium in (\ref{eq:stop_criterion})
selects the model with the best trade-off, according to $\delta$, between data adherence and complexity. The latter is defined as the number of factors.

As regards the choice of $\delta$, our solution hinges on the following scale-invariance property of the Kullback-Leibler divergence (see \cite{8447502}):

\begin{propo}
Let $x(t) \sim \mathcal{N}(0,\Sigma)$ , $t=1, \dots, N$ be i.i.d. \textit{random vectors} taking values in $\mathbb{R}^m$ and 
let $\mathrm{x}(t)$ be a realization of $x(t)$. 
Define the sample covariance estimator as 
\textcolor{black}{$\mathbf{ \hat{\Sigma}} := \frac{1}{N} \sum_{t=1}^{N} \mathrm{x}(t)\mathrm{x}^\top(t).
$} The Kullback-Leibler divergence between $\Sigma$ and $\mathbf{\hat{\Sigma}}$ is a random variable whose distribution depends only on the number $N$ of random variables and on the dimension $m$ of each random variable. Namely
\textcolor{black}{
\begin{equation}
    d := \mathcal{D}_{KL} (\Sigma \Vert \mathbf{ \hat{\Sigma} }) =   \frac{1}{2} (\log |Q_N| + \tr (Q_{N}^{-1}) - m)
\end{equation} 
}
where $Q_N$ is the random matrix defined by $Q_N:= \frac{1}{N} \sum_{t=1}^{N} {\tilde{x}(t) \tilde{x} (t)^\top }$ with $\tilde{x}(t)$ being i.i.d. normalized Gaussian random vectors: $\tilde{x}(t) \sim \mathcal{N}(0, I_m).$
\end{propo}

In view of this result, we can empirically approximate the distribution of the random variable $ d = \mathcal{D}_{KL} (\Sigma \Vert \mathbf{\hat{\Sigma}}_{r}) $  by a standard Monte Carlo method. 
 In particular, after choosing a probability $\alpha \in (0, 1) $ and $N$, we can find the neighborhood of radius $\delta_{\alpha}$ (in the Kullback-Leibler topology)  for which $ \textit{Pr} ( d \leq \delta_{\alpha} ) = \alpha $ .

\section{Numerical simulations} 
\label{sec:Simul}
To provide empirical evidence of the estimation performance of the algorithm, simulations studies have been performed by using the software Matlab-R2019b. 

\begin{figure*}
  \centering
\begin{tabular}{c c c}
    \includegraphics[width=0.4\textwidth, height = 0.16\textheight]{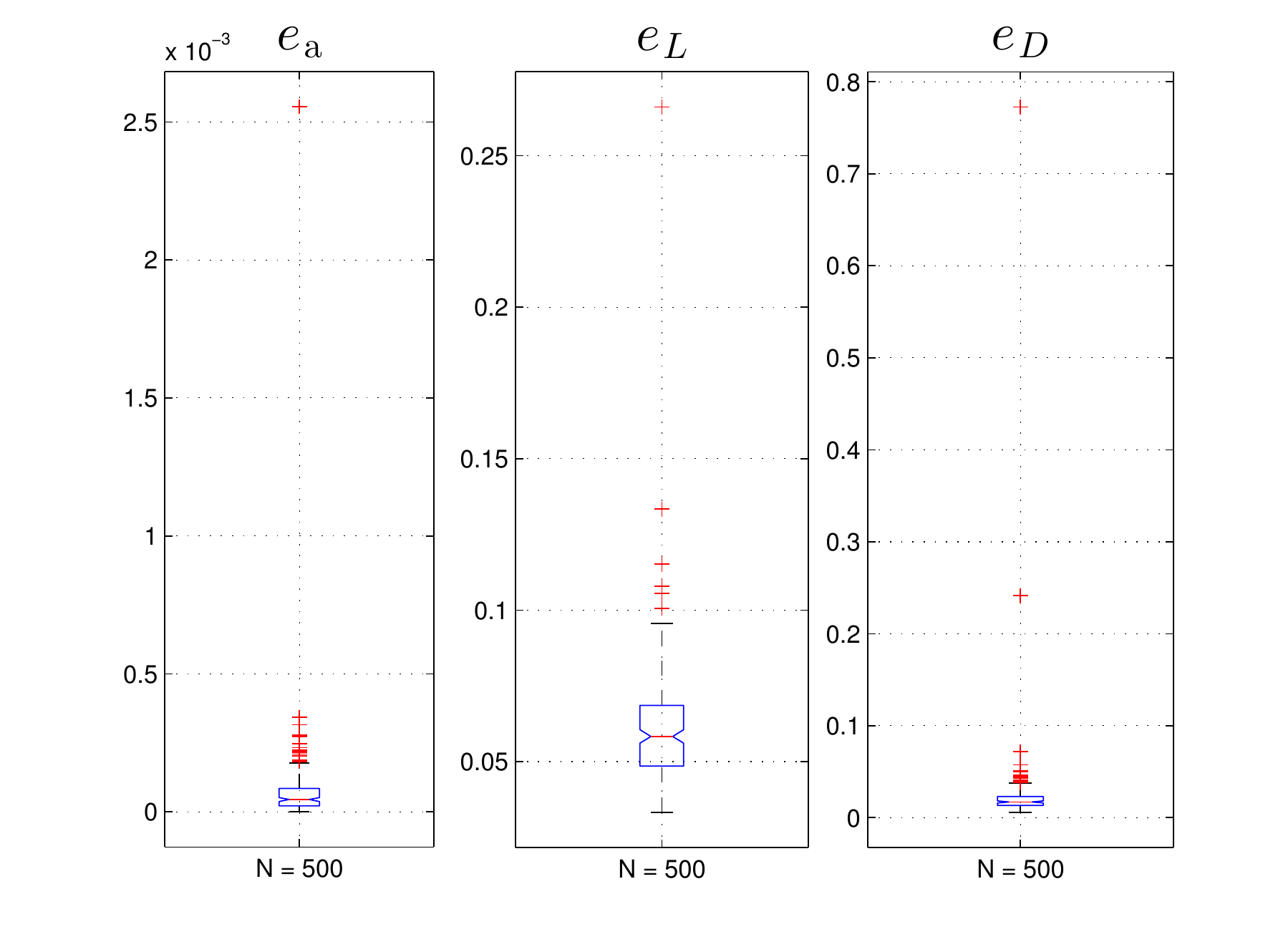} &
    \includegraphics[width=0.4\textwidth, height = 0.16\textheight]{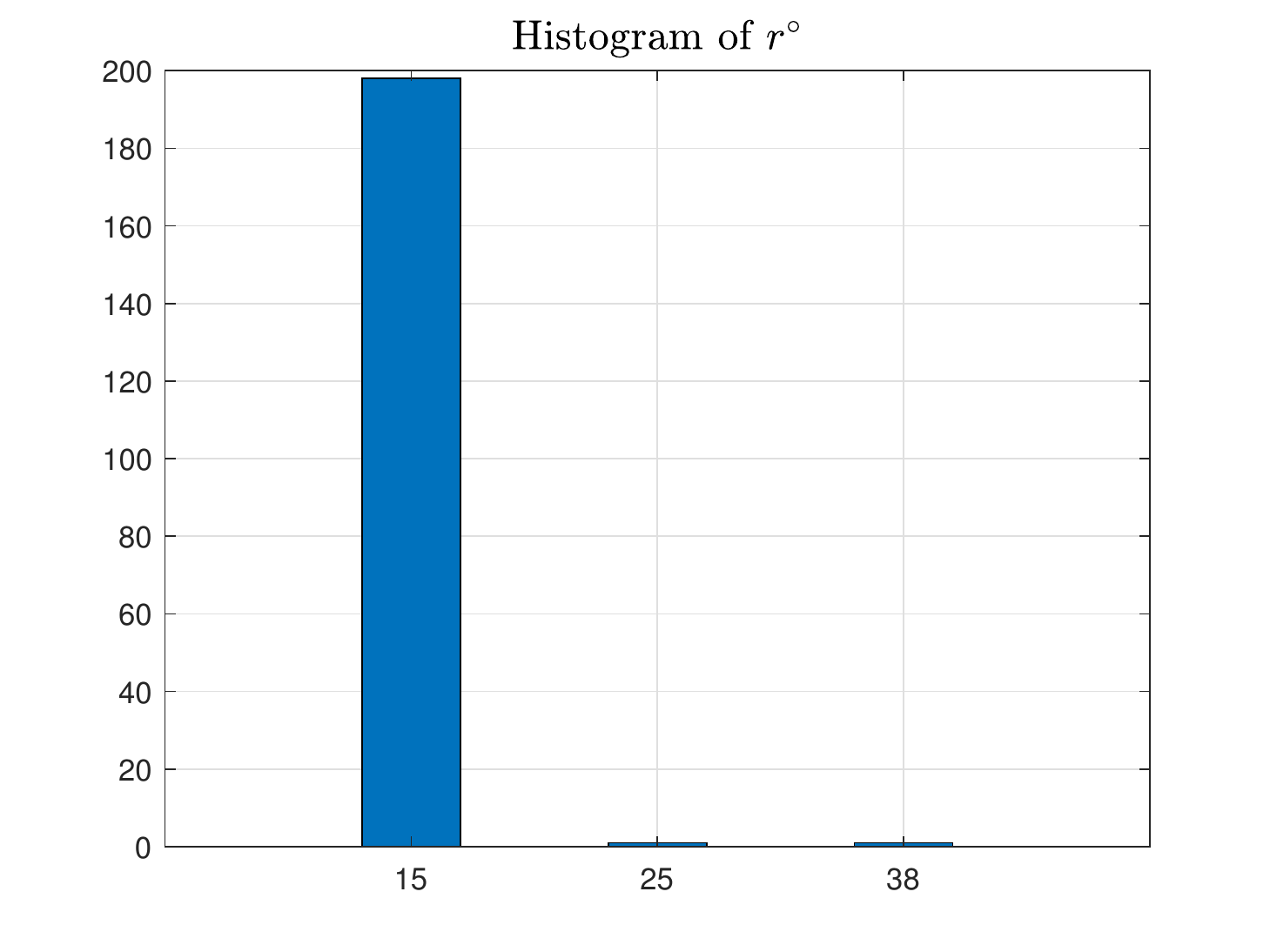} \\
   (a) & (b) 
 \end{tabular}
\caption{ Case $m=100$, $r=15$, $p=4$ with $N=500$. (a)
Box-plot of the errors $e_{\mathrm{a}}$, $e_{\mathrm{L}}$ and $e_{\mathrm{D}}$; (b) bar-plot of estimated $r^\circ$.}
  \label{fig:train_step}
\end{figure*}

We considered the case of a covariance matrix $\Sigma$, computed as the sum of a randomly generated positive semi-definite low-rank matrix $L$ of dimension $m$ and rank $r$, and a randomly generated positive
definite diagonal matrix $D$ such that $\|D\|_F\approx\| L\|_F$, i.e. the idiosyncratic noise is not negligible. Furthermore, we generated $a(z)$ by randomly choosing $p$ stable poles; without loss of generality we fixed $a_0=1$. Regarding the parameters of our procedure, we set $\alpha=0.99$, $\varepsilon=0.03$, $l_{max}= 200$, $\varepsilon_s=10^{-6}$ and $i_{max}=200$. In what follows, we analyze the following quantities:
\begin{itemize}
    \item the relative error on $\mathrm{a}$, $ e_{\mathrm{a}} := \Vert  \mathrm{a}   - \mathrm{a}^\circ_{r^\circ} \Vert / \Vert  \mathrm{a} \Vert$.
     \vspace{0.5mm}
    \item the relative error on $L$, $e_L := \Vert   L - L^\circ_{ r^\circ} \Vert_{F} / \Vert    L \Vert_{F}$;
    \vspace{0.5mm}
    \item the relative error on $D$, $e_D := \Vert    D -  {D}^\circ_{r^\circ} \Vert_{F} / \Vert    D \Vert_{F}$;
    \vspace{2mm}
  \end{itemize}

{\em First study:} We performed $200$ Monte Carlo runs with $m = 40$, $r = 10$ and \textcolor{black}{$p=5$}. Figure \ref{fig:m40_errTau}, \ref{fig:m40_errL} and \ref{fig:m40_errD} show that the proposed algorithm reaches good results since the errors are on the order of $10^{-1}$ for $L$ and $D$, and on the order of $10^{-3}$ for $\mathrm{a}$. The original low-rank and diagonal matrices are recovered with negligible numerical errors. The bar-plot of Figure \ref{fig:m40_bar} shows good performances also in the estimation of $r$. The worst case is achieved when $N = 200$: the efficiency is just the $53\%$ since in the $44\%$ of the cases the rank is underestimated as $ r^\circ  = 9$. On the other hand,  with $N = 500$ and $N = 800$ the efficiency reaches the $93\%$ and $95.5\%$ respectively, with just some outliers.

{\em Second study:} We performed $200$ Monte Carlo runs with $m = 100$, $r = 15$ and $p=4$. Hence, we considered a high dimensional case. The errors are plotted in Figure 2a and in Figure 2b the bar-plot of the estimated rank.
The estimates are very good. In particular the rank in the majority of the cases is correctly estimated. A similar study with essentially the same results has been conducted with the same parameters except for the order $p=5$ of the AR dynamics.

\textcolor{black}{It is remarkable that in all the simulations the algorithm stops before reaching the maximum number of iterations $l_{max}$  when $r^o$ approximates quite well the true value, as in these cases condition \eqref{eq:stop_condition} is satisfied.}

\section{Conclusions} \label{sec:Conclusions} In this paper we have considered the problem to estimate an AR factor model. More precisely, we have proposed an approximate moment matching procedure which alternates a static factor analysis step and an AR identification step. Empirical results showed that the algorithm estimates accurately the number of factors. This paradigm can be generalized to the case in which the order of the AR process is unknown. In such a scenario one could choose a criterium with complexity term, such as BIC. It is clear this requires to  compare the candidate models over a two dimensional grid (one dimension is $r$ and the other one is $p$), as a consequence the computational burden will be increased.


\end{document}